\documentclass[11pt]{article}
\usepackage{color}

\usepackage{amsmath,amssymb,mleftright,mathtools}
\usepackage{graphicx}
\usepackage{theorem}
\usepackage{here}
\usepackage{comment}

\newtheorem{thm}{Theorem}[section]

\newtheorem{cor}{Corollary}[section]
\newtheorem{prp}{Proposition}[section]
\theorembodyfont{\rmfamily}

\makeatletter

\@addtoreset{equation}{section}
\makeatother
\def\al{{\alpha}}
\def\be{{\beta}}

\def\de{{\delta}}

\def\si{{\sigma}}

\def\psih{{\widehat \psi}}

\def\Th{{\Theta}}
\def\De{{\Delta}}
\def\Si{{\Sigma}}

\def\Om{{\Omega}}

\def\La{{\Lambda}}

\def\Deh{{\widehat \De}}
\def\bTh{{\text{\boldmath $\Th$}}}

\def\bSi{{\text{\boldmath $\Si$}}}

\def\bOm{{\text{\boldmath $\Om$}}}
\def\bLa{{\text{\boldmath $\La$}}}
\def\bXi{{\text{\boldmath $\Xi$}}}

\def\bPsi{{\text{\boldmath $\Psi$}}}
\def\bPhi{{\text{\boldmath $\Phi$}}}
\def\bThh{{\widehat \bTh}}
\def\bSih{{\widehat \bSi}}

\def\bSit{{\widetilde \bSi}}

\def\e{{\text{\boldmath $e$}}}

\def\F{{\text{\boldmath $F$}}}

\def\H{{\text{\boldmath $H$}}}
\def\I{{\text{\boldmath $I$}}}

\def\O{{\text{\boldmath $O$}}}

\def\Q{{\text{\boldmath $Q$}}}
\def\R{{\text{\boldmath $R$}}}
\def\S{{\text{\boldmath $S$}}}

\def\X{{\text{\boldmath $X$}}}
\def\Y{{\text{\boldmath $Y$}}}

\def\Ybt{{\widetilde \Y}}

\def\Nc{{\cal N}}
\def\Oc{{\cal O}}

\def\Wc{{\cal W}}

\def\tr{{\rm tr\,}}
\def\diag{{\rm diag\,}}

\def\[{{\text{\boldmath $[$}}}
\def\]{{\text{\boldmath $]$}}}

\def\zero{{\bf\text{\boldmath $0$}}}

\def\|{{\,|\,}}

\def\/{{\Bigr/\!\!}}

\def\1r{{\rm (1)}}
\def\2r{{\rm (2)}}
\def\3r{{\rm (3)}}
\def\4r{{\rm (4)}}
\def\5r{{\rm (5)}}

\def\non{{\nonumber}}
\def\bXit{{\widetilde \bXi}}

\begin{document}
\title{Generalized Bayes Estimators with Closed forms for the Normal Mean and Covariance Matrices}
\author{Ryota Yuasa\footnote{Graduate School of Economics, University of Tokyo, 
7-3-1, Hongo, Bunkyo-ku, Tokyo 113-0033, Japan \quad
E-Mail: ryooyry@yahoo.co.jp}
\quad and \quad
Tatsuya Kubokawa\footnote{Faculty of Economics, University of Tokyo, 
7-3-1, Hongo, Bunkyo-ku, Tokyo 113-0033, Japan \quad
E-Mail: tatsuya@e.u-tokyo.ac.jp }}
\maketitle
\begin{abstract}
In the estimation of the mean matrix in a multivariate normal distribution, the generalized Bayes estimators with closed forms are provided, and the sufficient conditions for their minimaxity are derived relative to both matrix and scalar quadratic loss functions.
The generalized Bayes estimators of the covariance matrix are also given with closed forms, and the dominance properties are discussed for the Stein loss function.

\par\vspace{4mm}
\end{abstract}

\section{Introduction}
\label{sec:int}

The problems of estimating the mean and covariance matrices in the multivariate normal distribution  are addressed in a decision-theoretic framework.
These correspond to the estimation of the regression-coefficients matrix and the covariance matrix of the error terms in multivariate linear regression models.
Efron and Morris (1972) provided the empirical Bayes estimator of the mean matrix and showed the minimaxity relative to the quadratic loss.
Since then, dominance properties of shrinkage estimators have been studies in the literature.
For the good account of this problem, see Fourdrinier, Strawderman and Wells (2018), Tsukuma and Kubokawa (2020).
In the decision-theoretic framework, an ultimate goal is to establish the minimaxity of generalized Bayes estimators, because it is close to being admissible.
When the covariance matrix is known, Tsukuma (2009), Matsuda and Komaki (2015), and Tsukuma and Kubokawa (2017) showed the minimaxity of the generalized Bayes estimators.
In the case of unknown covariance matrix, however, no dominance properties have established.
Although Matsuda and Komaki (2015) and Tsukuma and Kubokawa (2017) established the minimaxity by checking the super-harmonicity of the prior distributions based on Stein (1973, 81) for known variance, the same method cannot be applied in the case of unknown covariance matrices.
It is also difficult to evaluate the risk functions directly, because the estimators involve integrals with respect to a symmetric matrix.

\medskip
In this paper, we derive the generalized Bayes shrinkage estimators of the mean matrix with closed forms against the specific prior distributions.
This approach is an extension of Maruyama and Strawderman (2005) who gave closed-forms of the generalized Bayes estimators in estimation of a mean vector.
Since the generalized Bayes estimators of the mean matrix do not contain integrals, it is easy to evaluate their risk functions, and we can derive conditions for their minimaxity.
The same specific prior distributions can also produce the generalized Bayes estimators of the covariance matrix with closed forms. 

\medskip
More specifically, we consider the canonical form of multivariate linear regression model: $\X$ and $\S$ are $m\times p$ and $p\times p$ random matrices, respectively, such that $\X$ and $\S$ are mutually independently distributed as $\X\sim \Nc_{m\times p}(\bTh, \I_m\otimes \bSi)$ and $\S \sim \Wc_p(n, \bSi)$ for $n\geq p$.
Then we derive the generalized Bayes estimators of $\bTh$ and $\bSi$ and show that for a specific prior, they are simplified as 
\begin{align*}
{\widehat \bTh}^{GB}&=
\X - k_0 \X\{\I_p+(1-k_0)\S^{-1}\X^\top\X\}^{-1},\\
\bSih^{GB}&=\bSit-k_0\bSit\{\I_p+(1-k_0)\S^{-1}\X^\top\X\}^{-1},
\end{align*}
for $\bSit=(m+c+1)^{-1}\S$ and $k_0$ given in (\ref{al}).
These derivations are given in Sections \ref{sec:mean} and \ref{sec:cov}.

\medskip
It is interesting to examine whether the generalized Bayes shrinkage estimators improve on the unbiased estimators which are used as benchmark estimators.
To this end, one first derives unbiased estimators of the risk functions using the Stein identity for the normal distribution and the Stein-Haff identity for the Wishart distribution.
We can employ this approach and find conditions for dominance of the generalized Bayes estimators, because they are given in closed forms.
In Section \ref{sec:minimax}, we derive conditions for ${\widehat \Th}^{GB}$ to dominate the unbiased estimator $\X$, where the estimators are evaluated relative to the two loss functions of the matrix quadratic loss and the scalar quadratic loss.
The dominance properties under the matrix quadratic loss have been studied in Bilodeau and Kariya (1989), Abu-Shanab, Kent and Strawderman (2012) and Matsuda and Strawderman (2021).

\medskip
Concerning the estimation of the covariance matrix $\bSi$, in Section \ref{sec:minimax}, we discuss the dominance properties of $\bSih^{GB}$ relative to the Stein loss function.
For $p>m$, we use the approach based on the unbiased risk estimation to get the condition for $\bSih^{GB}$ to dominate the unbiased estimator.
In the case of $m\geq p$, however, we cannot show the dominance property of $\bSih^{GB}$ based on the unbiased risk estimation, because the unbiased risk estimate of $\bSih^{GB}$ is larger than the risk of the unbiased estimator at some extreme cases.
We also consider the simultaneous estimation of $(\bTh, \bSi)$ relative to the Kullback-Leibler divergence, and the conditions for the dominance of $(\bThh^{GB}, \bSih^{GB})$ can be provided.
Especially, in the case of  $m\geq p$, we can borrow the risk gain of  $\bThh^{GB}$ to show the improvement of $(\bThh^{GB}, \bSih^{GB})$.

\medskip
Some numerical investigations based on simulation experiments are given in Section \ref{sec:sim} and the dominance results of the generalized Bayes estimators are supported numerically.
The paper is concluded with some remarks in Section \ref{sec:remark}.

\section{Generalized Bayes Estimators of the Mean Matrix with Closed forms}
\label{sec:mean}

Let $\X$ and $\S$ be $m\times p$ and $p\times p$ random matrices, respectively, such that $\X$ and $\S$ are mutually independently distributed as $\X\sim \Nc_{m\times p}(\bTh, \I_m\otimes \bSi)$ and $\S \sim \Wc_p(n, \bSi)$ for $n\geq p$.
We consider the estimation of the mean matrix $\bTh$ based on $\X$ and $\S$.
As loss functions, we treat the following two losses:
\begin{align}
L_M(\bThh, \bTh)=&(\bThh-\bTh)\bSi^{-1}(\bThh-\bTh)^{\top},
\label{mloss}\\
L_Q(\bThh, \bTh)=&\tr\{(\bThh-\bTh)\bSi^{-1}(\bThh-\bTh)^{\top}\},
\label{sloss}
\end{align}
which are called the matrix quadratic and scalar quadratic loss functions, respectively.
The estimator is evaluated by the corresponding risk functions.

\medskip
When a prior distribution is assumed for $(\bTh, \bSi)$, the (generalized) Bayes estimator relative to the two losses is expressed in the same form as 
$$
\bThh^B = E[\bTh\bSi^{-1}\mid \X,\S] \{ E[\bSi^{-1}\mid \X,\S]\}^{-1},
$$
where $E[\cdot\mid \X,\S]$ is the expectation with respect to the posterior distribution of $(\bTh, \bSi)$ given $(\X, \S)$.
We want to obtain generalized Bayes estimators of $\bTh$ which improves on $\X$.
To this end, we assume the following hierarchical prior distribution which was suggested in Section 4 of Tsukuma (2009):
\begin{equation}
\begin{split}
\bTh \mid \bSi, \bOm &\sim \Nc_{m\times p}(\zero, \bOm \otimes \bSi),\\
\pi(\bOm) &\propto |\I_m+\bOm|^{-a/2-m} |\bOm|^{c/2},\\
\pi(\bSi^{-1}) &\propto |\bSi^{-1}|^{(b-1)/2}.
\end{split}
\label{prior1}
\end{equation}

We first consider the case of $p>m$.
Let $\R$ be an $m\times m$ orthogonal matrix such that 
\begin{equation}
\X\S^{-1}\X^\top=\R\F\R^\top\quad {\rm and}\quad \F=\diag(f_1, \ldots, f_m),
\label{dec1}
\end{equation}
with $f_1>\cdots>f_m>0$.
Throughout the paper, we define constant $k_0$ by
\begin{equation}
k_0={a-c+p+(m\wedge p)-1 \over a+m+p+(m\wedge p)},
\label{al}
\end{equation}
which is less than one.

\begin{thm}
\label{thm:BE1}
In the case of $p>m$, assume that $-2<c<a+p$ and $b>-n-m-1$.
The generalized Bayes estimator against the prior $(\ref{prior1})$ is
\begin{equation}
\bThh^{GB1}=\X - \R\F^{-1}E[\bLa\mid\F]\{\F^{-1}(\I_m+\F)-E[\bLa\mid\F]\}^{-1}\R^\top\X,
\label{GB1}
\end{equation}
where $E[\bLa\mid\F]$ is the expectation of $\bLa$ with respect to the density
\begin{align}
\pi(\bLa\mid \F)\propto |\bLa|^{(a-c+m+p-1)/2-(m+1)/2}&|\I_m-\bLa|^{(c+m+1)/2-(m+1)/2} \non\\&|\I_m-\F(\I_m+\F)^{-1}\bLa|^{(b-a-m+n)/2}
I(0<\bLa<\I_m),
\label{posterior1}
\end{align}
for the indicator function $I(0<\bLa<\I_m)$.

When $b=a+m-n$, assume that $-2<c<a+p$, $a>-2m-1$.
Then the generalized Bayes estimator is expressed in a closed form as
\begin{align}
\bThh^{GB}&=\X - k_0 \{\I_m+(1-k_0)\X\S^{-1}\X^\top\}^{-1}\X
\non\\&=\X - k_0 \X\{\I_p+(1-k_0)\S^{-1}\X^\top\X\}^{-1}.
\label{GB1c}
\end{align}
\end{thm}

We next consider the case of $m\geq p$.
Let $\Q$ be a $p\times p$ nonsingular matrix such that 
\begin{equation}
\Q^\top\S\Q=\I_p, \quad \Q^\top\X^\top\X\Q=\F \quad {\rm and}\quad \F=\diag(f_1, \ldots, f_p),
\label{dec2}
\end{equation}
with $f_1>\cdots>f_p>0$.

\begin{thm}
\label{thm:BE2}
In the case of $m\geq p$, assume that $p-m-2<c<a+p$, $b>-n-m-1$.
The generalized Bayes estimator against the prior $(\ref{prior1})$ is
\begin{equation}
\bThh^{GB2}=\X - \X\Q\F^{-1}E[\bLa\mid\F]\{\F^{-1}(\I_p+\F)-E[\bLa\mid\F]\}^{-1}\Q^{-1},
\label{GB2}
\end{equation}
where $E[\bLa\mid\F]$ is the expectation of $\bLa$ with respect to the density
\begin{align}
\pi(\bLa\mid \F)\propto |\bLa|^{(a-c+2p-1)/2-(p+1)/2}&|\I_p-\bLa|^{(c+m+1)/2-(p+1)/2} \\&|\I_p-\F(\I_p+\F)^{-1}\bLa|^{(b-a-p+n)/2}
I(0<\bLa<\I_p).
\label{posterior2}
\end{align}

When $b=a+p-n$, assume that $p-m-2<c<a+p$.
Then the generalized Bayes estimator is expressed in a closed form as
\begin{align}
\bThh^{GB}&=\X - k_0 \X\{\I_p+(1-k_0)\S^{-1}\X^\top\X\}^{-1}.
\label{GB2c}
\end{align}
\end{thm}

It is noted that the generalized Bayes estimators (\ref{GB1c}) and (\ref{GB2c}) are of the same form in the two cases.
The proofs of Theorems \ref{thm:BE1} and \ref{thm:BE2} are given below.

\medskip
{\it Proof of Theorem \ref{thm:BE1}}.\ \ 
Since the conditional distribution of $\bTh$ given $\bSi$, $\bOm$, $\X$ and $\S$ is
$$
\bTh\mid \bSi, \bOm, \X, \S \sim \Nc_{m\times p}((\I_m+\bOm^{-1})^{-1}\X, (\I_m+\bOm^{-1})^{-1}\otimes \bSi) ,
$$
the generalized Bayes estimator is $\bThh^{GB}=\X - E[(\I_m+\bOm)^{-1}\X\bSi^{-1}\mid \X, \S]\{E[\bSi^{-1}\mid\X,\S]\}^{-1}$, where
\begin{align}
\pi(\bSi^{-1}, \bOm\mid\X,\S)\propto&
|\bSi^{-1}|^{(b+m+n+p)/2-(p+1)/2}{ |\bOm|^{c/2}\over |\I_m+\bOm|^{a/2+m+p/2}}\non\\
&\times \exp\big[-2^{-1}\tr\{ \bSi^{-1}(\S+\X^\top(\I_m+\bOm)^{-1}\X)\}\big],
\label{post0}
\end{align}
for $b+n+m>-1$.
After integrating the conditional expectations with respect to $\bSi^{-1}$, the estimator is written as
\begin{align*}
\bThh^{GB}=\X - &E[(\I_m+\bOm)^{-1}\X\{\S+\X^\top(\I_m+\bOm)^{-1}\X\}^{-1}\mid \X, \S]\\
&\times \big\{E[\{\S+\X^\top(\I_m+\bOm)^{-1}\X\}^{-1}\mid\X,\S]\big\}^{-1},
\end{align*}
where
$$
\pi(\bOm\mid\X,\S)\propto
{ |\bOm|^{c/2}\over |\I_m+\bOm|^{a/2+m+p/2}}{1\over |\S+\X^\top(\I_m+\bOm)^{-1}\X|^{(b+m+n+p)/2}}.
$$
Making the transformation $\bXi=(\I_m+\bOm)^{-1}$ with $d\bOm=|\bXi|^{-(m+1)}d\bXi$, we rewrite  the posterior density and $\bThh^{GB}$ as
\begin{equation}
\pi(\bXi\mid\X,\S)\propto
 |\S+\X^\top\bXi\X|^{-(b+m+n+p)/2} |\bXi|^{(a-c+p)/2-1}|\I_m-\bXi|^{c/2},
\label{postp}
\end{equation}
and 
\begin{align}
\bThh^{GB}=&\X - E[\bXi\X\{\S+\X^\top\bXi\X\}^{-1}\mid \X, \S]\{E[\{\S+\X^\top\bXi\X\}^{-1}\mid\X,\S]\}^{-1}\non\\
=& \X- E[\bXi\X\S^{-1/2}\{\I_p+\S^{-1/2}\X^\top\bXi\X\S^{-1/2}\}^{-1}\mid \X, \S]\non\\
&\quad\quad\times\big\{E[\{\I_p+\S^{-1/2}\X^\top\bXi\X\S^{-1/2}\}^{-1}\mid\X,\S]\big\}^{-1}\S^{1/2}.
\label{GBp}
\end{align}

We now consider the case of $p>m$.
Then, we show that $\bThh^{GB}$ is expressed as $\bThh^{GB1}$ in (\ref{GB1}).
Let $\Oc(m)$ be a set of $m\times m$ orthogonal matrices.
Since , by the singular-value decomposition, there are $\O_m\in \Oc(m)$ and $\O_p\in \Oc(p)$ such that $\X\S^{-1/2}=\O_m [\F^{1/2}, \zero_{m\times(p-m)}]\O_p$.
Then, it is demonstrated that 
\begin{align*}
|\I_p+\S^{-1/2}\X^\top\bXi\X\S^{-1/2}|=&|\I_p + \O_p^\top[\F^{1/2}, \zero_{m\times(p-m)}]^\top\O_m^\top \bXi \O_m [\F^{1/2}, \zero_{m\times(p-m)}]\O_p|\\
=&|\I_m + \F \O_m^\top \bXi \O_m|,
\end{align*}
and 
\begin{align*}
E[&\bXi\X\S^{-1/2}\{\I_p+\S^{-1/2}\X^\top\bXi\X\S^{-1/2}\}^{-1}\mid \X, \S]
\big\{E[\{\I_p+\S^{-1/2}\X^\top\bXi\X\S^{-1/2}\}^{-1}\mid\X,\S]\big\}^{-1}\S^{1/2}
\\
=&
E[\bXi\O_m [\F^{1/2}, \zero_{m\times(p-m)}]\O_p\{\I_p + \O_p^\top[\F^{1/2}, \zero_{m\times(p-m)}]^\top\O_m^\top \bXi \O_m [\F^{1/2}, \zero_{m\times(p-m)}]\O_p\}^{-1}\mid \X, \S]\\
&\times\big\{E[\{\I_p + \O_p^\top[\F^{1/2}, \zero_{m\times(p-m)}]^\top\O_m^\top \bXi \O_m [\F^{1/2}, \zero_{m\times(p-m)}]\O_p\}^{-1}\mid\X,\S]\big\}^{-1}\S^{1/2}\\
=&
\O_m E[\O_m^\top\bXi\O_m \{\F^{-1} + \O_m^\top \bXi \O_m\}^{-1}\mid \X, \S]\\
&\quad\times\big\{E[\{\F^{-1} + \O_m^\top \bXi \O_m \}^{-1}\mid\X,\S]\big\}^{-1}[\F^{1/2}, \zero_{m\times(p-m)}]\O_p\S^{1/2}.
\end{align*}
Note that $\O_m=\R$ and $[\F^{1/2}, \zero_{m\times(p-m)}]\O_p\S^{1/2}=\O_m^\top \X=\R^\top \X$.
Making the transformation $\bXit=\O_m^\top \bXi \O_m$ gives
\begin{align*}
&\pi(\bXit\mid\X,\S)\propto
 |\F^{-1}+\bXit|^{-(b+m+n+p)/2} |\bXit|^{(a-c+p)/2-1}|\I_m-\bXit|^{c/2},\\
&\bThh^{GB}=
\X- \R E[\bXit \{\F^{-1} + \bXit\}^{-1}\mid \X, \S]\big\{E[\{\F^{-1} + \bXit \}^{-1}\mid\X,\S]\big\}^{-1}\R^\top\X.
\end{align*}
Making the transformation $\bLa=(\I_m+\F)^{1/2}(\F+\bXit^{-1})^{-1}(\I_m+\F)^{1/2}$, we can see that $0<\bLa<\I_m$ and $d\bLa=|\I_m+\F|^{(m+1)/2}|\F\bXi+\I_m|^{-(m+1)}d\bXi$.
Then, $\pi(\bXi\mid\F)d\bXi=\pi(\bLa\mid\F)d\bLa$, where
\begin{align}
\pi(\bLa\mid \F)\propto &|\bLa|^{(a-c+m+p-1)/2-(m+1)/2}|\I_m-\bLa|^{(c+m+1)/2-(m+1)/2} \non\\
&\times|\I_m-\F(\I_m+\F)^{-1}\bLa|^{(-a+b-m+n)/2}
I(0<\bLa<\I_m),
\label{Mp1}
\end{align}
for $a-c+p>0$ and $c>-2$.
Also, we have
\begin{align*}
E[&\bXit \{\F^{-1} + \bXit\}^{-1}\mid \X, \S]\big\{E[\{\F^{-1} + \bXit \}^{-1}\mid\X,\S]\big\}^{-1}
\\
=&
(\I_m+\F)^{-1/2}E[\bLa\mid\F](\F^{-1}(\I_m+\F)-E[\bLa\mid\F])^{-1}(\I_m+\F)^{1/2}\F^{-1},
\end{align*}
which leads to the expression given in (\ref{GB1}) if $E[\bLa\mid\F]$ is diagonal matrix. Now we can show that $\int_{0<\bLa<\I_p}\{\bLa\times |\bLa|^{(a-c+m+p-1)/2-(m+1)/2}|\I_m-\bLa|^{(c+m+1)/2-(m+1)/2} |\I_m-\F(\I_m+\F)^{-1}\bLa|^{(a-b-m+n)/2}\} d\bLa$ is diagonal matrix. Making the transformation $\bLa^*=\H_i\bLa\H_i$, where $\H_i=\I_m-2\e_i\e_i^{\top}, \e_i$ is the $m$ dimensional fundamental unit vector, then we have
\begin{align*}
&\int_{0<\bLa<\I_m}\{\bLa\times |\bLa|^{(a-c+m+p-1)/2-(m+1)/2}|\I_m-\bLa|^{(c+m+1)/2-(m+1)/2}  \\ &\hspace{3cm}|\I_m-\F(\I_m+\F)^{-1}\bLa|^{(a-b-m+n)/2}\} d\bLa\\=&\int_{0<\bLa^*<\I_m}\{\H_i\bLa^*\H_i\times |\bLa^*|^{(a-c+m+p-1)/2-(m+1)/2}|\I_m-\bLa^*|^{(c+m+1)/2-(m+1)/2} \\ &\hspace{3cm}|\I_m-\F(\I_m+\F)^{-1}\bLa^*|^{(a-b-m+n)/2}\} d\bLa^*.
\end{align*}
By comparing the right hand side and left hand side, we see that the off-diagonal elements of the $i$-th column and $i$-th vector are zero. Therefore $E[\bLa\mid\F]$ is diagonal.

In the case of $b=a+m-n$, the posterior distribution (\ref{posterior1}) is the multivariate beta distribution $Beta_m(a-c+m+p-1, c+m+1)$, so that we have $E[\bLa\mid\F]=k_0\I_m$ for $k_0=(a-c+m+p-1)/\{(a-c+m+p-1)+(c+m+1)\}=(a-c+m+p-1)/(a+2m+p)$.
Hence, the estimator (\ref{GB1}) is simplified as the first equality in  (\ref{GB1c}).
Note that
\begin{equation}\label{eq1}
(\I_m+\be\X\S^{-1}\X^\top)^{-1}=\I_m - \be\X(\S+\be\X\X^\top)^{-1}\X^\top
\end{equation}
for constant $\be$.
Using the equality (\ref{eq1}), we can further rewrite the estimator (\ref{GB1c}) as
\begin{align*}
&\X-k_0\{\I_m+(1-k_0)\X\S^{-1}\X^\top\}^{-1}\X
=\X-k_0\Big\{\I_m-\X\Big({1\over1-k_0}\S+\X^\top\X\Big)^{-1}\X^\top\Big\}\X\\
&=\X-k_0\X\Big({1\over1-k_0}\S+\X^\top\X\Big)^{-1}{1\over 1-k_0}\S
=\X-k_0\X\big\{\I_p+(1-k_0)\S^{-1}\X^\top\X\big\}^{-1},
\end{align*}
which leads to the expression given in the second equality of (\ref{GB1c}).
\hfill$\Box$

\medskip
{\it Proof of Theorem \ref{thm:BE2}}.\ \ 
We next consider the case of $m \geq p$.
We begin with the posterior density (\ref{postp}) and the generalized Bayes estimator (\ref{GBp}).
Then, we show that $\bThh^{GB}$ is expressed as $\bThh^{GB2}$ in (\ref{GB2}).
By the singular-value decomposition, there are $\O_m\in \Oc(m)$ and $\O_p\in \Oc(p)$ such that $\X\S^{-1/2}=\O_m [\F^{1/2}, \zero_{p\times (m-p)}]^\top\O_p$.
Then, it is demonstrated that 
\begin{align*}
|\I_p+\S^{-1/2}\X^\top\bXi\X\S^{-1/2}|=&|\I_p + \O_p^\top[\F^{1/2}, \zero_{p\times (m-p)}]\O_m^\top \bXi \O_m [\F^{1/2}, \zero_{p\times (m-p)}]^\top\O_p|\\
=&|\I_p + \F (\O_m^\top \bXi \O_m)_{11}|,
\end{align*}
where $(\O_m^\top \bXi \O_m)_{11}$ is a $p\times p$ matrix in the decomposition
\begin{equation}
\O_m^\top \bXi \O_m = \begin{pmatrix}(\O_m^\top \bXi \O_m)_{11} & (\O_m^\top \bXi \O_m)_{12}\\ (\O_m^\top \bXi \O_m)_{21}&(\O_m^\top \bXi \O_m)_{22}\end{pmatrix}.
\label{decomp}
\end{equation}
Also, it is seen that
\begin{align*}
E[&\bXi\X\S^{-1/2}\{\I_p+\S^{-1/2}\X^\top\bXi\X\S^{-1/2}\}^{-1}\mid \X, \S]\{E[\{\I_p+\S^{-1/2}\X^\top\bXi\X\S^{-1/2}\}^{-1}\mid\X,\S]\}^{-1}\S^{1/2}
\\
=&
E[\bXi\O_m [\F^{1/2}, \zero_{p\times (m-p)}]^\top\O_p\{\I_p + \O_p^\top[\F^{1/2}, \zero_{p\times (m-p)}]\O_m^\top \bXi \O_m [\F^{1/2}, \zero_{p\times (m-p)}]^\top\O_p\}^{-1}\mid \X, \S]\\
&\times\big\{E[\{\I_p + \O_p^\top[\F^{1/2}, \zero_{p\times (m-p)}]\O_m^\top \bXi \O_m [\F^{1/2}, \zero_{p\times (m-p)}]^\top\O_p\}^{-1}\mid\X,\S]\big\}^{-1}\S^{1/2}\\
=&
\O_m E[\O_m^\top\bXi\O_m [\I_p, \zero_{p\times (m-p)}]^\top \{\F^{-1} + (\O_m^\top \bXi \O_m)_{11}\}^{-1}\mid \X, \S]\\
&\quad\times\big\{E[\{\F^{-1} + (\O_m^\top \bXi \O_m)_{11} \}^{-1}\mid\X,\S]\big\}^{-1}\F^{1/2}\O_p\S^{1/2}.
\end{align*}
Note that $\Q=\S^{-1/2}\O_p^\top$, $\X=\O_m[\F^{1/2}, \zero_{p\times (m-p)}]^\top\Q^{-1}$.
Making the transformation $\bXit=\O_m^\top \bXi \O_m$ gives
\begin{align*}
&\pi(\bXit\mid\X,\S)\propto
 |\F^{-1}+\bXit_{11}|^{-(b+m+n+p)/2} |\bXit|^{(a-c+p)/2-1}|\I_m-\bXit|^{c/2},\\
&\bThh^{GB}=
\X- \O_m E[ [\bXit_{11}, \bXit_{21}^\top]^\top \{\F^{-1} +\bXit_{11}\}^{-1}\mid \X, \S]\big\{E[\{\F^{-1} + \bXit_{11} \}^{-1}\mid\X,\S]\big\}^{-1}\F^{1/2}\Q^{-1}.
\end{align*}
It is noted that $|\bXit|=|\bXit_{11}||\bXit_{22.1}|$ and $|\I_m-\bXit|=|\I_p-\bXit_{11}||\I_{m-p}-\bXit_{22.1}||\I_{m-p}-(\I_{m-p}-\bXit_{22.1})^{-1/2}\bXit_{21}(\bXit_{11}-\bXit_{11}^2)^{-1}\bXit_{21}^\top(\I_{m-p}-\bXit_{22.1})^{-1/2}|$ for $\bXit_{22.1}=\bXit_{22}-\bXit_{21}\bXit_{11}^{-1}\bXit_{21}^\top$, where $\bXit_{ij}$ is decomposed similarly to (\ref{decomp}) for $i, j = 1, 2$.
Writing $\pi(\bXit\mid\X,\S)d\bXit=\pi(\bXit_{11}, \bXit_{21}, \bXit_{22.1}\mid\X,\S)d\bXit_{11}d\bXit_{21}d\bXit_{22.1}$ where $\I_p>\bXit_{11}>\O_p, \I_{m-p}>\bXit_{22.1}>\O_{m-p}, \I_{m-p}-(\I_{m-p}-\bXit_{22.1})^{-1/2}\bXit_{21}\{\bXit_{11}(\I_p-\bXit_{11})\}^{-1}\bXit_{21}^{\top}(\I_{m-p}-\bXit_{22.1})^{-1/2}>\O_{m-p}$, we can see that $\pi(\bXit_{11}, -\bXit_{21}, \bXit_{22.1}\mid\X,\S)=\pi(\bXit_{11}, \bXit_{21}, \bXit_{22.1}\mid\X,\S)$, which implies that
$$
E[ [\bXit_{11}, \bXit_{21}^\top]^\top \{\F^{-1} +\bXit_{11}\}^{-1}\mid \X, \S]
= E[ [\bXit_{11}, \zero_{(m-p)\times m}^\top]^\top \{\F^{-1} +\bXit_{11}\}^{-1}\mid \X, \S].
$$
Thus, the generalized Bayes estimator and the related posterior density can be rewritten as
\begin{align*}
&\bThh^{GB}=
\X- \X\Q \F^{-1/2}E[ \bXit_{11}\{\F^{-1} +\bXit_{11}\}^{-1}\mid \X, \S]\big\{E[\{\F^{-1} + \bXit_{11} \}^{-1}\mid\X,\S]\big\}^{-1}\F^{1/2}\Q^{-1},\\
&\pi(\bXit_{11}\mid\X,\S)\propto
 |\F^{-1}+\bXit_{11}|^{-(b+m+n+p)/2} |\bXit_{11}|^{(a-c+m)/2-1}|\I_p-\bXit_{11}|^{(c+m-p)/2}.
\end{align*}
Making the transformation $\bLa=(\I_p+\F)^{1/2}(\F+\bXit_{11}^{-1})^{-1}(\I_p+\F)^{1/2}$, we can see that $0<\bLa<\I_p$ and $d\bLa=|\I_p+\F|^{(p+1)/2}|\F\bXit_{11}+\I_p|^{-(p+1)}d\bXit_{11}$.
Then, $\pi(\bXit_{11}\mid\F)d\bXit_{11}=\pi(\bLa\mid\F)d\bLa$, where
\begin{align}
\pi(\bLa\mid \F)\propto &|\bLa|^{(a-c+2p-1)/2-(p+1)/2}|\I_p-\bLa|^{(c+m+1)/2-(p+1)/2} \non\\
&\times |\I_p-\F(\I_p+\F)^{-1}\bLa|^{(b-a-p+n)/2}
I(0<\bLa<\I_p),
\label{Mp2}
\end{align}
for $a-c+p>0$ and $c+m>p-2$.
Since $E[\bXit_{11} \{\F^{-1} + \bXit_{11}\}^{-1}\mid \X, \S]\big\{E[\{\F^{-1} + \bXit_{11} \}^{-1}\mid\X,\S]\big\}^{-1}=(\I_p+\F)^{-1/2}E[\bLa\mid\F](\F^{-1}(\I_p+\F)-E[\bLa\mid\F])^{-1}(\I_p+\F)^{1/2}\F^{-1}$, one gets the expression given in (\ref{GB2}) if $E[\bLa\mid\F]$ is diagonal matrix, which can be verified by using the same arguments as in the proof of Theorem \ref{thm:BE1}. 
In the case of $b=a+p-n$, similarly to Theorem \ref{thm:BE1}, we have $E[\bLa\mid\F]=k_0\I_p$ for $k_0=(a-c+2p-1)/(a+m+2p)$, so that the estimator (\ref{GB2}) is simplified as in (\ref{GB2c}).
\hfill$\Box$

\section{Generalized Bayes Estimators of Covariance Matrix with Closed forms}
\label{sec:cov}

In this section, we derive the generalized Bayes estimators of the covariance matrix with closed forms.
The problem we consider is the estimation of $\bSi$ in the model described in the beginning of Section \ref{sec:mean}, namely, $\X\sim \Nc_{m\times p}(\bTh, \I_m\otimes \bSi)$ and $\S \sim \Wc_p(n, \bSi)$ for $n\geq p$.
When estimator $\bSih$ is evaluated by the risk function relative to the Stein loss function
\begin{equation}
L_S(\bSih, \bSi) = \tr(\bSih\bSi^{-1}) - \log |\bSih\bSi^{-1}| - p,
\label{steinloss}
\end{equation}
the (generalized) Bayes estimator of $\bSi$ against a prior distribution of $(\bTh, \bSi)$ is given by $\bSih^B=(E[\bSi^{-1}|\X, \S])^{-1}$.
We assume the same prior distribution as in (\ref{prior1}), and derive the generalized Bayes estimators in the two cases of $p>m$ and $m\geq p$.

\begin{thm}
\label{thm:CBE}
{\rm (i)}\ In the case of $p>m$, assume that $-2<c<a+p$ and $b>-n-m-1$.
The generalized Bayes estimator against the prior $(\ref{prior1})$ is
\begin{equation}
\bSih^{GB1}={1\over b+m+n+p}\left[ \S + \X^\top\R\left\{ (\I_m+\F)(E[\bLa\mid\F])^{-1}-\F\right\}^{-1}\R^\top\X\right],
\label{CGB1}
\end{equation}
where $\R$ and $\F$ are defined in $(\ref{dec1})$, and $E[\bLa\mid\F]$ is the expectation of $\bLa$ with respect to the density in $(\ref{posterior1})$.
When $b=a+m-n$, assume that $-2<c<a+p$, $a>-2m-1$.
Then the generalized Bayes estimator is expressed in a closed form as
\begin{align}
\bSih^{GB}&={1\over a+2m+p}\left[ \S + k_0 \X^\top\left\{\I_m +(1-k_0)\X\S^{-1}\X^\top\right\}^{-1}\X\right]
\non\\&=\bSit-k_0\bSit(\I_p+(1-k_0)\S^{-1}\X^\top\X)^{-1}
\label{CGB1c}
\end{align}
where $\bSit=(m+c+1)^{-1}\S$ and $k_0$ is given in $(\ref{al})$.

\smallskip
{\rm (ii)}\ In the case of $m\geq p$, assume that $p-m-2<c<a+p$, $b>-n-m-1$.
The generalized Bayes estimator against the prior $(\ref{prior1})$ is
\begin{equation}
\bSih^{GB2}={1\over b+m+n+p}\left[ \S + (\Q^\top)^{-1}\left\{ (\F^{-1}+\I_p)(E[\bLa\mid\F])^{-1}-\I_p\right\}^{-1}\Q^{-1}\right],
\label{CGB2}
\end{equation}
where $\Q$ and $\F$ are defined in $(\ref{dec2})$, and $E[\bLa\mid\F]$ is the expectation of $\bLa$ with respect to the density in $(\ref{posterior2})$.
When $b=a+p-n$, assume that $p-m-2<c<a+p$.
Then the generalized Bayes estimator is expressed in a closed form as
\begin{align}
\bSih^{GB}&={1\over a+m+2p}\left[ \S + k_0 \left\{(\X^\top\X)^{-1}+(1-k_0)\S^{-1}\right\}^{-1}\right]\non
\\&=\bSit-k_0\bSit(\I_p+(1-k_0)\S^{-1}\X^\top\X)^{-1}
\label{CGB2c}
\end{align}
\end{thm}

{\it Proof}.\ \ From the posterior distribution of $\bSi^{-1}$ given in (\ref{post0}), the generalized Bayes estimator of $\bSi$ is
\begin{align}
\bSih^{GB}=&{1\over b+m+n+p}\big[ E[\{\S+\X^\top(\I_m+\bOm)^{-1}\X\}^{-1}\mid \X, \S]\big]^{-1} \non\\
=&
{1\over b+m+n+p}\S^{1/2}\big[ E[\{\I_p+\S^{-1/2}\X^\top\bXi\X\S^{-1/2}\}^{-1}\mid \X,\S]\big]^{-1}\S^{1/2},
\label{Sip1}
\end{align}
where the posterior density of $\bXi=(\I_m+\bOm)^{-1}$ is given in (\ref{postp}).

\smallskip
For part (i), the same arguments as in the proof of Theorem \ref{thm:BE1} are used. 
Since $\X\S^{-1/2}$ is decomposed as $\X\S^{-1/2}=\O_m [\F^{1/2}, \zero_{m\times(p-m)}]\O_p$ for $\O_m\in \Oc(m)$ and $\O_p\in \Oc(p)$, we can rewrite $\bSih^{GB}$ as
\begin{equation}
\bSih^{GB}
=
{1\over b+m+n+p}\S^{1/2}\O_p^\top\begin{pmatrix} \F^{1/2}\{E[(\F^{-1}+\bXit)^{-1}\mid \F]\}^{-1}\F^{1/2} & \zero\\ \zero&\I_{p-m}\end{pmatrix}\O_p\S^{1/2},
\label{Sip2}
\end{equation}
for $\bXit=\O_m^\top\bXi\O_m$.
Making the transformation $\bLa=(\I_m+\F)^{1/2}(\F+\bXit^{-1})^{-1}(\I_m+\F)^{1/2}$ and noting that $E[\bLa\mid\F]$ is a diagonal matrix, we have
\begin{align*}
\F^{1/2}&\{E[(\F^{-1}+\bXit)^{-1}\mid \F]\}^{-1}\F^{1/2} \non\\
=&\F^{1/2}\{\F-\F(\I_m+\F)^{-1/2}E[\bLa\mid \F](\I_m+\F)^{-1/2}\F\}^{-1}\F^{1/2}
\non\\
=&\F^{1/2}\{\F-\F(\I_m+\F)^{-1/2}E[\bLa\mid \F](\I_m+\F)^{-1/2}\F\}^{-1}\F^{1/2}
\\
=&\I_m + \F^{1/2}(\I_m+\F-\F E[\bLa\mid\F])^{-1}E[\bLa\mid\F]\F^{1/2},
\end{align*}
where the posterior density of $\bLa$ is given in (\ref{Mp1}).
Then from (\ref{Sip2}), it follows that
\begin{align*}
\bSih^{GB}
=&
{1\over b+m+n+p}\Big[\S + \S^{1/2}\O_p^\top\begin{pmatrix}\F^{1/2}\\ \zero\end{pmatrix}\O_m^\top \O_m (\I_m+\F-\F E[\bLa\mid\F])^{-1}\\
&\hspace{4cm}\times E[\bLa\mid\F]\O_m^\top \O_m (\F^{1/2}, \zero)\O_p\S^{1/2}\Big]
\non\\
=&
{1\over b+m+n+p}\left\{\S + \X^\top \O_m (\I_m+\F-\F E[\bLa\mid\F])^{-1}E[\bLa\mid\F]\O_m^\top \X\right\}.
\end{align*}
Noting that $\O_m=\R$, one gets the first expression in (\ref{CGB1}).
The second expression in  (\ref{CGB1}) can be verified by using the equality (\ref{eq1}).

\smallskip
For part (ii), we can use the same arguments as in the proof of Theorem \ref{thm:BE2}. 
Since $\X\S^{-1/2}$ is decomposed as $\X\S^{-1/2}=\O_m [\F^{1/2}, \zero_{p\times(m-p)}]^\top\O_p$ for $\O_m\in \Oc(m)$ and $\O_p\in \Oc(p)$, we can rewrite $\bSih^{GB}$ in (\ref{Sip1}) as
\begin{equation}
\bSih^{GB}
=
{1\over b+m+n+p}\S^{1/2} \left\{ E[( \I_p + \O_p^\top \F^{1/2}\bXit_{11}\F^{1/2}\O_p)^{-1} \mid \F] \right\}^{-1}\S^{1/2},
\label{Sip3}
\end{equation}
for $\bXit=\O_m^\top\bXi\O_m$.
Making the transformation $\bLa=(\I_p+\F)^{1/2}(\F+\bXit_{11}^{-1})^{-1}(\I_p+\F)^{1/2}$ and noting that $E[\bLa\mid\F]$ is a diagonal matrix, we have
\begin{align*}
\big\{ E[&( \I_p + \O_p^\top \F^{1/2}\bXit_{11}\F^{1/2}\O_p)^{-1} \mid \F] \big\}^{-1}\\
=& \O_p^\top\{ \I_p - (\I_p+\F)^{-1}\F E[\bLa\mid \F]\}^{-1}\O_p\\
=& \I_p + \O_p^\top (\F^{-1}+\I_p - E[\bLa\mid \F])^{-1} E[\bLa\mid \F]\O_p,
\end{align*}
where the posterior density of $\bLa$ is given in (\ref{Mp2}).
Noting that $\S^{1/2}\O_p^\top=\S \S^{-1/2}\O_p=(\Q^\top)^{-1}\Q^{-1}\Q=(\Q^\top)^{-1}$, from (\ref{Sip3}), we can see that
\begin{align*}
\bSih^{GB}
=&
{1\over b+m+n+p}\Big[\S + \S^{1/2}\O_p^\top(\F^{-1}+\I_p - E[\bLa\mid \F])^{-1} E[\bLa\mid \F]\O_p\S^{1/2}\Big]
\non\\
=&
{1\over b+m+n+p}\left\{\S + (\Q^\top)^{-1} (\F^{-1}+\I_p - E[\bLa\mid \F])^{-1} E[\bLa\mid \F]\Q^{-1}\right\},
\end{align*}
which leads to the first expression in (\ref{CGB2}).
The second expression in (\ref{CGB2}) can be easily derived.
\hfill$\Box$

\section{Dominance Properties}
\label{sec:minimax}

\subsection{Estimation of the mean matrix}
We now show that the generalized Bayes estimator with the closed form improve on $\X$ relative to the matrix and scalar quadratic loss functions.
To this end, we begin by providing the unbiased estimators of the risk functions of the general class of shrinkage estimators
\begin{equation}
\bThh^{SH}=\left\{
\begin{array}{ll}
(\I_m-\R\bPhi(\F)\R^\top)\X, & {\rm if}\ p>m,\\
\X(\I_p-\Q\bPhi(\F)\Q^{-1}), & {\rm if}\ m\geq p,
\end{array}\right.
\label{SH}
\end{equation}
where $\bPhi(\F)=\diag(\phi_1(\F), \ldots, \phi_{\ell}(\F))$ and $\F=\diag(f_1, \ldots, f_\ell)$ for $\ell=m\wedge p$.
In the case of $p>m$, we can measure the risk of the estimator $\bThh^{SH}$ with respect to the matrix quadratic loss $L_M(\bThh, \bTh)$ in (\ref{mloss}).
Tsukuma and Kubokawa (2015, 20) derived the unbiased estimator of the risk function given by ${\widehat R}_M(\bThh^{SH})=p\I_m + \R\bPhi^*\R^\top$ for $\bPhi^*=\diag(\phi_1^*, \ldots, \phi_m^*)$, where 
\begin{align}
\phi_i^*=& (n-p+2\ell-3)f_i\phi_i^2-2(p-m+1)\phi_i -4f_i^2\phi_i{\partial\phi_i\over \partial f_i}- 4f_i{\partial\phi_i\over \partial f_i}
\non\\
&-2\sum_{j\not= i}^m{f_i^2\phi_i^2\over f_i-f_j}+2\sum_{j\not= i}^m{f_i\phi_i f_j\phi_j\over f_i-f_j}-2\sum_{j\not= i}^m{f_i\phi_i- f_j\phi_j\over f_i-f_j}.
\label{uerm}
\end{align}
Concerning the scalar quadratic loss $L_Q(\bThh, \bTh)$ in (\ref{sloss}), we can provide the unified expression of the unbiased estimator of the risk function in the cases of $p>m$ and $m\geq p$, which is given by ${\widehat R}_S(\bThh^{SH})=mp + \tr(\bPhi^*)$, where
\begin{align}
\tr(\bPhi^*)=& \sum_{i=1}^\ell \Big\{(n-p+2\ell-3)f_i\phi_i^2-2(|p-m|+1)\phi_i -4f_i^2\phi_i{\partial\phi_i\over \partial f_i}- 4f_i{\partial\phi_i\over \partial f_i}\Big\}
\non\\
&-2 \sum_{i=1}^\ell \sum_{j=i+1}^\ell\Big( {f_i^2\phi_i^2-f_j^2\phi_j^2\over f_i-f_j}+2{f_i\phi_i- f_j\phi_j\over f_i-f_j}\Big).
\label{uers}
\end{align}
For the details of the derivation, see Konno (1990, 91, 92), Tsukuma (2009) and Tsukuma and Kubokawa (2020).

\medskip
Using (\ref{uerm}) and (\ref{uers}), we derive conditions for the generalized Bayes estimators in (\ref{GB1c}) and (\ref{GB2c}) to improve on $\X$.
More generally, we consider the class of estimators
\begin{equation}
{\widehat \bTh}^{G} =
\X - \al \X\{\I_p+\be\S^{-1}\X^\top\X\}^{-1},
\label{GTh}
\end{equation}
which corresponds to $\phi_i^G=\al/(1+\be f_i)$, where $\al$ and $\be$ are positive constants.

\begin{thm}
\label{thm:minimax}
{\rm (i)}\ \ Assume that $p>m+1$ and $n\geq p$.
If $\al$ and $\be$ satisfy the condition 
\begin{equation}
\al \leq 2(p-m-1)\be/(n-p+2m+1),
\label{Gmin1}
\end{equation}
then the estimator $\bThh^G$ dominates $\X$ relative to the matrix quadratic loss function $L_M(\bThh, \bTh)$ in $(\ref{mloss})$. 

{\rm (ii)}\ \ Assume that $|p-m|>1$ and $n\geq p$.
If $\al$ and $\be$ satisfy the condition 
\begin{equation}
\al \leq 2(|p-m|-1)\be/(n-p+2\ell+1),
\label{Gmin2}
\end{equation}
then the estimator $\bThh^G$ dominates $\X$ relative to the scalar quadratic loss function $L_Q(\bThh, \bTh)$ in $(\ref{sloss})$. 
\end{thm}

It is noted that the generalized Bayes estimator $\bThh^{GB}$ given in (\ref{GB1c}) or (\ref{GB2c}) corresponds to the case of $\al=k_0$ and $\be=1-k_0$ for $k_0$ given in (\ref{al}).
The conditions for the dominance of $\bThh^{GB}$ follows from Theorem \ref{thm:minimax}.

\begin{cor}
\label{cor:minimax1}
{\rm (i)}\ \ For the prior distribution $(\ref{prior1})$, assume that $p>m+1$, $n\geq p$, $b=a+m-n$, $a>-2m-1$, $-2<c<a+p$ and 
\begin{equation}
k_0={a-c+m+p-1\over a+2m+p} \leq {2(p-m-1)\over n+p-1}.
\label{minc1}
\end{equation}
Then the generalized Bayes estimator $(\ref{GB1c})$ improves on $\X$ relative to the matrix quadratic loss function $L_M(\bThh, \bTh)$ in $(\ref{mloss})$. 

{\rm (ii)}\ \ For the prior distribution $(\ref{prior1})$, assume that $n\geq p$, $b=a+(m\wedge p)-n$, $a>-m-(m\wedge p)-1$, $m\wedge p-m-2<c<a+p$ and
\begin{equation}
k_0={a-c+p+(m\wedge p)-1\over a+m+p+(m\wedge p)} \leq {2(|p-m|-1)\over n-p+2(m\wedge p)+2|p-m|-1}.
\label{minc2}
\end{equation}
Then the generalized Bayes estimators $(\ref{GB1c})$ and $(\ref{GB2c})$ improve on $\X$ relative to the scalar quadratic loss function $L_Q(\bThh, \bTh)$ in $(\ref{sloss})$. 
\end{cor}

{\it Proof of Theorem \ref{thm:minimax}}.\ \ 
For part (i), let $A=n-p+2m-3$, $B=p-m+1$ and $g_i=1/\{1+\be f_i\}$ for simplicity.
Then $\phi_i^{G}=\al g_i$ and $\partial \phi_i^{G}/\partial f_i = -\al\be g_i^2$.
The former part in (\ref{uerm}) can be written as
\begin{align}
R_{M1, i}^{G}=&Af_i(\phi_i^{G})^2 - 2 B\phi_i^{G} - 4 f_i^2\phi_i^{G}{\partial \phi_i^{G}\over \partial f_i}- 4f_i{\partial \phi_i^{G}\over \partial f_i}\non\\
=& \al^2Af_ig_i^2-2\al Bg_i+4\al^2\be f_i^2g_i^3+4\al\be f_ig_i^2
\non\\
=& \al g_i\big\{ 4\al\be(f_ig_i)^2 + (\al A+4\be)f_ig_i -2B\big\}.
\label{minp1}
\end{align}
The latter part in (\ref{uerm}) can be written as
\begin{align}
R_{M2, i}^{G}=&-2\sum_{j\not= i}{1\over f_i-f_j}\big\{ f_i^2 (\phi_i^{G})^2-f_if_j\phi_i^{G}\phi_j^{G}+f_i\phi_i^{G}-f_j\phi_j^{G}\big\}
\non\\
=&-2\sum_{j\not= i}{\al\over f_i-f_j}\big\{ \al f_i^2 g_i^2-\al f_if_jg_i g_j+f_ig_i-f_jg_j\big\}\non\\
=& -2 \sum_{j\not= i}\al g_i^2g_j\{(\al+\be)f_i+1\}.
\label{minp2}
\end{align}
Note that $\phi_i^*$ in (\ref{uerm}) is $R_{M1,i}^{G}+R_{M2,i}^{G}$ and $R_{M2, i}^{G}\leq 0$.
Since $f_ig_i\leq 1/\be$, it is observed that
\begin{align}
R_{M1, i}^{G}=& \al g_i\big\{ 4\al\be(f_ig_i)^2 + (\al A+4\be)f_ig_i -2B\big\}\non\\
\leq& \al g_i \Big\{ 4{\al\over \be} + (\al A + 4\be){1\over \be} -2B\Big\}\non\\
=& {\al\over \be} g_i \{ (A+4)\al - 2(B-2)\be\}.
\label{minp3}
\end{align}
Thus, $R_{M1,i}^{G}\leq 0$ holds if $(A+4)\al\leq 2(B-2)\be$.

\medskip
For part (ii), we show that $\sum_{i=1}^\ell (R_{M1,i}^{G}+R_{M2,i}^{G})\leq 0$, where $B$ is replaced with $A=n-p+2\ell-3$ and $B=|p-m|+1$.
Using the same arguments as in the part (i), we can see that the dominance result holds if $(A+4)\al\leq 2(B-2)\be$.
\hfill$\Box$

\subsection{Estimation of the covariance matrix}
We next investigate a dominance property of the generalized Bayes estimators with the closed forms of the covariance matrix $\bSi$ relative to the Stein loss (\ref{steinloss}).
To this end, we begin by providing the unbiased estimators of the risk functions of the general class of estimators
\begin{equation}
\bSih^{SH}=\left\{
\begin{array}{ll}
\bSih_0 - n^{-1}\X^\top\R\F^{-1}\bPsi(\F)\R^\top\X, & {\rm if}\ p>m,\\
\bSih_0 - n^{-1}(\Q^\top)^{-1}\bPsi(\F)\Q^{-1}, & {\rm if}\ m\geq p,
\end{array}\right.
\label{SHC}
\end{equation}
where $\bSih_0=n^{-1}\S$, $\bPsi(\F)=\diag(\psi_1(\F), \ldots, \psi_{\ell}(\F))$ and $\F=\diag(f_1, \ldots, f_\ell)$ for $\ell=m\wedge p$.
The risk difference of the estimators $\bSih^{SH}$ and $\bSih_0$ is written as $\De_S(\bSih^{SH}, \bSih_0)=R_S(\bSih^{SH}, \bSi)-R_{S}(\bSih_0, \bSi)$, and Tsukuma and Kubokawa (2016) provided the unbiased estimator $\Deh_S(\bSih^{SH}, \bSih_0)$ of $\De_S(\bSih^{SH}, \bSih_0)$, where
\begin{align}
\Deh(\bSih^{SH}, \bSih_0)
={1\over n}\sum_{i=1}^\ell\Big\{ -d_i\psi_i+2f_i{\partial\psi_i\over \partial f_i} + 2 \sum_{j>i}^\ell{\psi_i-\psi_j\over f_i-f_j}f_j-n\log(1- \psi_i)\Big\},
\label{uerc}
\end{align}
where $d_i=n-p+2i-1$ for $i=1,\ldots,\ell$.
From Theorem \ref{thm:CBE}, the generalized Bayes estimators are expressed as 
$$
\bSih^{GB}=\left\{
\begin{array}{ll}
(a+2m+p)^{-1}\big[ \S + k_0 \X^\top\big\{\I_m +(1-k_0)\X\S^{-1}\X^\top\big\}^{-1}\X\big], & {\rm for}\ p>m,\\
(m+c+1)^{-1}\big[\S - k_0\S\{\I_p+(1-k_0)\S^{-1}\X^\top\X\}^{-1}, & {\rm for}\ m\geq p,
\end{array}
\right.
$$
which belong to the above class when $a+2m+p=n$ for $p>m$ and when $m+c+1=n$ for $m\geq p$.
We below show that the dominance properties are different in the two cases of $p>m$ and $m\geq p$.

\medskip
We first treat the case of $p>m$ and consider the general class of estimators
\begin{equation}
\bSih^{G1}=\bSih_0 + n^{-1}\al \X^\top\big\{\I_m +\be\X\S^{-1}\X^\top\big\}^{-1}\X,
\label{GSi1}
\end{equation}
where $\al$ and $\be$ are nonnegative constants, and $\bSih_0=n^{-1}\S$ is an unbiased estimator of $\bSi$.

\begin{thm}\label{thm:GSi1}
Assume that $p>m$.
The estimator $\bSih^{G1}$ dominates $\bSih_0$ relative to the Stein loss $(\ref{steinloss})$ if
\begin{equation}\label{GSi1c}
{\al\over\be}\leq{2(p-m)\over n-p+m}.
\end{equation}
\end{thm}

{\it Proof.}\ \ 
The estimator $\bSih^{G1}$ corresponds to $\psi_i^{G1}=-\al f_ig_i$ for $g_i=1/(1+\be f_i)$.
Note that 
$$
{\partial \psi_i^{G1}\over \partial f_i} \leq 0, \quad{\rm and}\quad
\sum_{j>i}^m {\psi_i^{G1}-\psi_j^{G1}\over f_i-f_j}f_j\leq 0.
$$
Using these facts, we can evaluate (\ref{uerc}) as
$$
\Deh(\bSih^{G1}, \bSih_0)
\leq{1\over n}\sum_{i=1}^m\{ \al d_i f_ig_i - n \log(1+\al f_ig_i)\}.
$$
It is here noted that  for $x_1\geq \cdots \geq x_k$, 
\begin{equation}
2\sum_{i=1}^k i x_i \leq (k+1)\sum_{i=1}^k x_i,
\label{ineq}
\end{equation}
which was provided in Dey and Srinivasan (1985).
Then, $\sum_{i=1}^m d_i f_ig_i\leq \sum_{i=1}^m (n-p+m)f_ig_i$, which is used to evaluate $\Deh(\bSih^{G1}, \bSih_0)$ from above as
$$
\Deh(\bSih^{G1}, \bSih_0)
\leq {1\over n}\sum_{i=1}^m\Big\{ (n-p+m)\al f_i g_i-n\log(1+ \al f_i g_i)\Big\}.
$$
Note that for $z>0$, $\log(1+z)\geq 2z/(2+z)$, which was used in Tsukuma and Kubokawa (2016).
Then, $-\log (1+\al/\be)\leq -2(\al/\be)/(2+\al/\be)$. 
It is also noted that $h(x)=(n-p+m)x-n\log(1+x)$ si convex in $x$.
Since $0\leq\al f_ig_i\leq\al/\be$, it holds that 
\begin{align}\label{SiG1p}
\Deh(\bSih^{G1}, \bSih_0)
&\leq \max\Big[0, \  {m\over n}\Big\{ (n-p+m){\al\over \be}-n\log\left(1+ {\al\over \be}\right)\Big\} \Big] \non \\
&\leq {m\over n}\max\Big[0,\  {\al/\be\over 2+\al/\be}\sum_{i=1}^m\Big\{ (n-p+m){\al\over \be} +2(-p+m)\Big\} \Big].
\end{align}
Therefore we can see that $\Deh\leq 0$ if ${\al/\be}\leq{2(p-m)/(n-p+m)}$, which leads to the condition (\ref{GSi1c}).
\hfill$\Box$

\bigskip
When $a+2m+p=n$, we can apply Theorem \ref{thm:GSi1} to find a condition for the generalized Bayes estimator to dominate the unbiased estimator, which corresponds to $\al=k_0$ and $\be=1-\al$.
For $a=n-2m-p$, the condition (\ref{GSi1c}) is 
$$
k_0={n-(m+c+1)\over n} \leq{2(p-m)\over n+p-m}
$$
for $n-(m+c+1)>0$.
Then the following corollary can be provided from Theorem \ref{thm:CBE}.

\begin{cor}
In the case of $p>m$, let $a+2m+p=n$.
The generalized Bayes estimator 
\begin{equation}\label{GBSi1}
\bSih^{GB}=\bSih_0 + n^{-1}\big[k_0 \X^\top\big\{\I_m +(1-k_0)\X\S^{-1}\X^\top\big\}^{-1}\X\big]
\end{equation}
dominates $\bSih_0$ relative to the Stein loss if $c$ satisfies that $\max\{-2, n(n-p+m)/(n+p-m)-m-1\}< c < n-2m$.
\end{cor}

We next consider the case of $m\geq p$. 
We begin by handling generalized Bayes estimator $\bSih^{GB}$ in (\ref{CGB2c}), which corresponds to $\psi_i^{GB}=1-c_0+c_0 k_0h_i$ for $c_0=n/(m+c+1)$ and  $h_i=1/\{1+(1-k_0)f_i\}$.
Then from (\ref{uerc}),
$$
\Deh(\bSih^{GB}, \bSih_0)
={1\over n}\sum_{i=1}^p\Big\{ -d_i\left(1-c_0+c_0k_0h_i\right) - 2c_0k_0(1-k_0) \sum_{j\geq i}^p f_j h_ih_j-n\log\left(c_0+c_0k_0h_i\right)\Big\}.
$$
To investigate a necessary condition for $\Deh(\bSih^{GB}, \bSih_0)\leq 0$, consider the extreme case of $f_i$'s, namely $f_i\to \infty$ for all $i$.
Since $h_i\to 0$ and $f_ih_i\to 1/(1-\al)$, it is clear that 
$$
\lim_{f_i's\to \infty} \Deh(\bSih^{GB}, \bSih_0)
=p( c_0 -\log c_0 -1 ),
$$
which is larger than or equal to zero.
This implies that we should take $c_0=1$, that is, $m+c+1=n$ is a necessary condition for $\Deh(\bSih^{GB}, \bSih_0)\leq 0$. 
Then $\psih_i^{GB}=k_0/\{1+(1-k_0)f_i\}$.

\medskip
More generally, we consider the estimator
\begin{equation}
\bSih^{G2}=\bSih_0 - \al\bSih_0(\I_p+\be\S^{-1}\X^\top\X)^{-1},
\label{GSi2}
\end{equation}
which corresponds to $\psih_i^{G2}=\al/(1+\be f_i)$ for nonnegative constatnts $\al$ and $\be$.
In this case, however, we cannot find any condition for the dominance of $\bSih^{G2}$ over $\bSih_0$.
From (\ref{uerc}),
\begin{equation}\label{GSip}
\Deh(\bSih^{G2}, \bSih_0)
={1\over n}\sum_{i=1}^p\Big\{ -d_i\al g_i-2\al\be f_i g_i^2 - 2\al\be \sum_{j>i}^p f_j g_ig_j-n\log(1- \al g_i)\Big\},
\end{equation}
for $g_i=1/(1+\be f_i)$.
To investigate a necessary condition for $\Deh(\bSih^{G2}, \bSih_0)\leq 0$, consider the extreme case of $f_i$'s, namely $f_i\to 0$ for all $i$.
Since $g_i\to 1$ as $f_i\to 0$, it is seen that
\begin{align*}
\lim_{f_i's\to 0} \Deh(\bSih^{G2}, \bSih_0)
=&{1\over n}\sum_{i=1}^p\{ -d_i\al +n\log(1- \al)\}
={p\over n}\{-n\al - n \log(1-\al)\}\\
=& p\{(1-\al) - \log(1-\al)-1\},
\end{align*}
which is positive for $\al>0$, because $(1-\al) - \log(1-\al)-1>0$ for $\al>0$.
This implies the following proposition.

\begin{prp}
\label{prp:Si}
Assume that $m\geq p$.
The does not exist any estimators in the class $\bSih^{G2}$ with $\al>0$ such that $\Deh(\bSih^{G2}, \bSih_0)\leq 0$ for any $f_i$'s.
Furthermore the generalized Bayes estimator $\bSih^{GB}$ does not satisfy the condition for $\Deh(\bSih^{GB}, \bSih_0)\leq 0$ for $f_i$'s close to 0.
\end{prp}

\subsection{Dominance result under the Kullback-Leibler divergence}
In the previous subsections, the generalized Bayes estimator of the mean matrix dominates the unbiased estimator, but such a dominance result cannot be guaranteed in the estimation of the covariance matrix.
These raise a question about whether any dominance result holds for the generalized Bayes estimators in simultaneous estimation of the mean and covariance matrices.
We here introduce the Kullback-Leibler divergence, which is given by
$$
D_{KL}(\phi(\Ybt|\bThh, \bSih) ,\phi(\Ybt|\bTh, \bSi))=\int \log \big\{{\phi(\Ybt|\bThh, \bSih)/ \phi(\Ybt|\bTh, \bSi)}\big\}\phi(\Ybt|\bThh, \bSih) d \Ybt,
$$
where $\phi(\Ybt|\bTh, \bSi)$ denotes the probability density function of $\Nc_{m\times p}(\bTh, \I_m\otimes \bSi)$.
A direct calculation shows that
\begin{align}
D_{KL}(\phi(\Ybt|\bThh, \bSih) ,\phi(\Ybt|\bTh, \bSi))
=&{m\over 2}\big\{\tr(\bSih\bSi^{-1})-\log(|\bSih\bSi^{-1}|)-p\big\}+{1\over 2}\tr\big\{(\bThh-\bTh)\bSi^{-1}(\bThh-\bTh)^\top\big\}\non\\
=&{m\over2}L_S(\bSih, \bSi)+{1\over 2}L_Q(\bThh, \bTh),\label{KL1}
\end{align}
and estimators $(\bThh, \bSih)$ are evaluated by the risk $E[D_{KL}(\phi(\Ybt|\bThh, \bSih) ,\phi(\Ybt|\bTh, \bSi))]$.

We begin by considering the general shrinkage estimators $(\bThh^{SH}, \bSih^{SH})$ given in (\ref{SH}) and (\ref{SHC}), and the risk difference between them and the unbiased estimators $(\X, \bSih_0)$ is denoted by $\De=\De((\bThh^{SH}, \bSih^{SH}), (\X, \bSih_0))=E[D_{KL}(\phi(\Ybt|\bThh^{SH}, \bSih^{SH}) ,\phi(\Ybt|\bTh, \bSi))]-E[D_{KL}(\phi(\Ybt|\X, \bSih_0) ,\phi(\Ybt|\bTh, \bSi))]$.
Then from (\ref{uers}) and (\ref{uerc}), the risk difference $\De$ is expressed as $\De=E[\Deh]$, where
\begin{align*}
2\Deh=& \sum_{i=1}^\ell \Big\{Af_i\phi_i^2-2B\phi_i -4f_i^2\phi_i{\partial\phi_i\over \partial f_i}- 4f_i{\partial\phi_i\over \partial f_i}\Big\}
\non\\
&-2 \sum_{i=1}^\ell \sum_{j=i+1}^\ell\Big( {f_i^2\phi_i^2-f_j^2\phi_j^2\over f_i-f_j}+2{f_i\phi_i- f_j\phi_j\over f_i-f_j}\Big)\\
&+{m\over n}\sum_{i=1}^\ell\Big\{ -d_i\psi_i+2f_i{\partial\psi_i\over \partial f_i} + 2 \sum_{j>i}^\ell{\psi_i-\psi_j\over f_i-f_j}f_j-n\log(1- \psi_i)\Big\},
\end{align*}
for $A=n-p+2\ell-3$, $B=|p-m|+1$ and $d_i=n-p+2i-1$.

\medskip
We first treat the case of $p>m$ and consider the estimators $({\widehat \bTh}^{G}, \bSih^{G1})$ given in (\ref{GTh}) and (\ref{GSi1}).

\begin{thm}
\label{thm:KLmin1}
Assume that $p>m+1$ and $n\geq p$.
If $\al$ and $\be$ satisfy the conditions 
\begin{equation}
{\al\over\be}\leq {2(p-m-1)\over n-p+2m+1},
\label{KLmin1}
\end{equation}
then the estimator $({\widehat \bTh}^{G}, \bSih^{G1})$  dominates $(\X, \bSih_0)$ relative to the Kullback-Leibler divergence.
\end{thm}

{\it Proof.}\ \ 
From (\ref{minp2}), ({\ref{minp3}) and (\ref{SiG1p}), the risk difference is written as
\begin{align}
2\Deh \leq& \sum_{i=1}^m \al g_i\Big\{ 4{\al\over\be} + {\al A+4\be\over \be} -2B\Big\}\non\\
&+ {m\over n}\max\left\{0, {\al/\be\over 2+\al/\be}\sum_{i=1}^m\Big\{ (n-p+m){\al\over \be} +2(-p+m)\Big\} \right\}.
\label{KL1p1}
\end{align}
The first term in RHS of (\ref{KL1p1}) is not positive if $\al/\be \leq 2(p-m-1)/(n-p+2m+1)$.
The second term in RHS of (\ref{KL1p1}) is not positive if $\al/\be\leq 2(p-m)/(n-p+m)$.
Since $2(p-m-1)/(n-p+2m+1)<2(p-m)/(n-p+m)$, it is sufficient to satisfy the condition $\al/\be\leq 2(p-m-1)/(n-p+2m+1)$.
\hfill$\Box$

\bigskip
Since the generalized Bayes estimator $({\widehat \bTh}^{GB}, \bSih^{GB})$ corresponds to the case of $\al=k_0$ and $\be=1-k_0$, Theorem \ref{thm:KLmin1} provides the condition for the dominance of the generalized Bayes estimator.

\begin{cor}
\label{cor:KLminGB1}
Assume that $p>m+1$.
Let $a=n-p-2m$ and $k_0=\{n-(c+m+1)\}/n$.
The estimator $({\widehat \bTh}^{GB}, \bSih^{GB})$  dominates $(\X, \bSih_0)$ relative to the Kullback-Leibler divergence if $c$ satisfies that $-2< c < n-2m$ and
\begin{equation}
c_{low}\coloneqq{n^2-(p-m)n-(p-1)(m+1)\over n+p-1}\leq c.
\label{KLminG1}
\end{equation}
There exists a $c$ satisfying these conditions if $n>(p-1)(m-1)/(2p-3m-1)$ and $p>(3m+1)/2$.
\end{cor}

\medskip
We next treat the case of $m\geq p$ and consider the estimators $({\widehat \bTh}^{G}, \bSih^{G2})$ given in (\ref{GTh}) and (\ref{GSi2}).
Although $\bSih^{G2}$ cannot improve on $\bSih_0$ as seen from Proposition \ref{prp:Si}, we can borrow the risk gain of ${\widehat \bTh}^{G}$ in the framework of simultaneous estimation of $(\bTh, \bSi)$ to establish the dominance property of $({\widehat \bTh}^{G}, \bSih^{G2})$.

\begin{thm}
\label{thm:KLmin}
Assume that $m\geq p$.
If $\al$ and $\be$ satisfy the inequalities $\al<1$ and
\begin{equation}
(n+p+1)\al + {m\al\be\over 2(1-\al)}\leq 2(m-p-1)\be+2{m\over n}\be,
\label{KLmin}
\end{equation}
then the estimator $({\widehat \bTh}^{G}, \bSih^{G2})$  dominates $(\X, \bSih_0)$ relative to the Kullback-Leibler divergence.
\end{thm}

{\it Proof.}\ \ 
From (\ref{minp1}), (\ref{minp2}), ({\ref{minp3}) and (\ref{GSip}), the risk difference is written as
\begin{align}
2\Deh =& \sum_{i=1}^p\Big[ \al g_i\big\{ 4\al\be(f_ig_i)^2 + (\al A+4\be)f_ig_i -2B\big\} - 2 \sum_{j\not= i}\al g_i^2g_j\{(\al+\be)f_i+1\}\Big]\non\\
&+ {m\over n}\sum_{i=1}^p \Big\{ -d_i\al g_i-2\al\be f_i g_i^2 - 2\al\be \sum_{j>i}^p f_j g_ig_j-n\log(1- \al g_i)\Big\}.\label{KLp1}
\end{align}
The inequalities in (\ref{ineq}) are used to evaluate (\ref{KLp1}). We also use the inequality $-\log(1-x)\leq x+x^2/\{2(1-x)\}$ for $0<x<1$, which is provided in Maruyama and Strawderman (2012).
Note that $\sum_{i=1}^p d_i g_i \geq \sum_{i=1}^p ng_i$ and $-\log (1-\al g_i)\leq \al g_i + \al^2g_i^2/\{2(1-\al g_i)\}\leq \al g_i + \al^2 g_i/\{2(1-\al)\}$.
Thus, the second term in (\ref{KLp1}) is evaluated as
\begin{align*}
\sum_{i=1}^p &\Big\{ -d_i\al g_i-2\al\be f_i g_i^2 - 2\al\be \sum_{j>i}^p f_j g_ig_j-n\log(1- \al g_i)\Big\}\\
&\leq   \sum_{i=1}^p \Big\{ -2\al\be f_i g_i^2 - 2\al\be \sum_{j>i}^p f_j g_ig_j + {n\al\over 2(1-\al)}\al g_i\Big\},
\end{align*}
so that 
\begin{align}
2\Deh \leq& \sum_{i=1}^p\Big[ \al g_i\big\{ 4\al\be(f_ig_i)^2 + (\al A+4\be)f_ig_i -2B\big\} - 2 \sum_{j\not= i}\al g_i^2g_j\{(\al+\be)f_i+1\}\Big]\non\\
&+ {m\over n}\sum_{i=1}^p \Big\{ -2\al\be f_i g_i^2 - 2\al\be \sum_{j>i}^p f_j g_ig_j + {n\al\over 2(1-\al)}\al g_i\Big\}.\label{KLp2}
\end{align}
Deleting several terms, we have
\begin{align}
2\Deh \leq& \sum_{i=1}^p \al g_i \Big\{ 4\al\be(f_ig_i)^2 + \Big(\al A+4\be-2 {m\over n}\be\Big)f_ig_i -2B  + {m\al\over 2(1-\al)}\Big\}.\label{KLp3}
\end{align}
Since $f_ig_i\leq 1/\be$, we can see that $\Deh\leq 0$ if
$$
{4\al\over \be} + \Big(\al A+4\be-2 {m\over n}\be\Big){1\over \be} -2B + {m\al\over 2(1-\al)}\leq 0,
$$
equivalently
$$
(A+4)\al + {m\al\be\over 2(1-\al)}\leq 2(B-2)\be+2{m\over n}\be,
$$
which leads to the condition (\ref{KLmin}).
\hfill$\Box$

\bigskip
Since the generalized Bayes estimator $({\widehat \bTh}^{GB}, \bSih^{GB})$ corresponds to the case of $\al=k_0$ and $\be=1-k_0$, Theorem \ref{thm:KLmin} provides the condition for the dominance of the generalized estimator.

\begin{cor}
\label{cor:KLminGB}
Assume that $m\geq p$.
Let $c=n-m-1$ and $k_0=(a+m+2p-n)/(a+m+2p)$.
The estimator $({\widehat \bTh}^{GB}, \bSih^{GB})$  dominates $(\X, \bSih_0)$ relative to the Kullback-Leibler divergence if $a$ satisfies that $n-m-p-1<a$ and
\begin{equation}
a\leq -m-2p+n\left(1+2{m-p-1+m/n \over n+m/2+p+1}\right) \eqqcolon a_{upp}.
\label{KLminG}
\end{equation}
There exists a $a$ satisfying these conditions if $m>(3p+1)/2$ and $n>(2p^2+mp-5m-1)/(2(2m-3p-1))$.
\end{cor}

\section{Simulation Study}
\label{sec:sim}

We now investigate the numerical performance of the suggested shrinkage estimators in the estimation of $\bTh$, $\bSi$ and $(\bTh, \bSi)$ and compare them with the benchmark estimators.
Among the suggested estimators, we treat the following three ones:

\smallskip
{\bf GB}: \ Concerning the generalized Bayes estimators $\bThh^{GB}$ and $\bSih^{GB}$, we choose parameters $a$ and $c$ as follows: \ $a=n-2m-p$ and $c=c_{low})$ when $p>m$ for $c_{low}$ given in $(\ref{KLminG1})$, and $a=a_{upp}$ and $c=n-m-1$ when $m\geq p$ for $a_{upp}$ given in $(\ref{KLminG})$. 
These estimators are denoted by GB.

\smallskip
{\bf G1}:\ Concerning the shrinkage estimators $\bThh^{G}$, $\bSih^{G1}$ and $\bSih^{G2}$ in $(\ref{GTh}), (\ref{GSi1}), (\ref{GSi2})$, we choose parameters $\al$ and $\be$ as follows: \ $\al=p-m-1$ and $\be=n-p+2m+1$ when $p>m+1$,  which satisfies the condition in $(\ref{Gmin2})$, and $\al=(m-p-1)/(n+m)$ and $\be=(n+p+1)/(n+m)$ when $m>p+1$,  which satisfies the condition in $(\ref{Gmin2})$ and $0<\al<1$. 
These estimators are denoted by G1.

\smallskip
{\bf G2}:\ In estimators $\bThh^{G}$, $\bSih^{G1}$ and $\bSih^{G2}$, another choice of parameters $\al$ and $\be$ is $\al=p-m, \be=n-p+m$ when $p>m$, which satisfies the condition in $(\ref{GSi1c})$, and $\al=(m-p)/(n+m-p), \be=n/(n+m-p)$ when $m\geq p$,  which satisfies $0\leq\al<1$.
These estimators are denoted by G2.

\smallskip
{\bf EM}:\ In this simulation, we add the Efron-Morris estimator to estimators we want to compare.
In the case of an unknown covariance matrix, as given in Tsukuma and Kubokawa (2015), the Efron-Morris estimator is $\bThh^{EM}=\X-c_{EM}\R\F^{-1}\R^\top\X$ where $c_{EM}=(|m-p|-1)/(\min\{n-p+2m, n+p\}+1)$ and $\R$ is an $m\times (m\wedge p)$ matrix with eigenvectors of $\X\S^{-1}\X^\top$. 
This estimator is denoted by EM.

\smallskip
The simulation experiments are conducted in the two cases of (1) $p=10, n=10, m=5, 25$ and (2) $p=10, n=25, m=5, 20, 30$, where the simulated data are generated from the multivariate normal distribution and the Wishart distribution with the following setup of parameters:
For singular values $s_1,\ldots,s_p$ of the mean matrix $\bTh$, we consider the case that $s_i=s_0+s_0\times (i-1)/(p/5-1)$ for $i\in\{1,\ldots, p/5\}$ and the others are $s_i=s_0/ 10^{q}$, where $s_0=0,1,10,20$ and $q=1, 1/2$, where the power $q$ controls the dispersion of eigenvalues, namely, eigenvalues for $q=1$ is more dispersed than that for $q=1/2$.
Eigenvectors of the mean matrix are constructed based on a matrix of which elements are generated from normal distribution with zero mean and one variance independently. 
For the covariance matrix $\bSi$, we consider the two cases of (1) $\si_{kl}=0.9+0.1\de_{kl}$ and (2) $\si_{kl}=0.5^{|k-l|}$, where $\de_{kl}=1$ when $k=l$ and $\de_{kl}=0$ when $k\neq l$.
It is noted that data are strongly correlated in the case (1). 
As averages based on the simulation experiments with 5,000 replications, we obtain the values of the percentage relative improvement in risk (PRIR) of $\bThh$ and/or  $\bSih$ over $\X$ and/or $\S/n$, where PRIRs in the three estimation problems are defined by 
\begin{align*}
&100\times {E[L_Q(\X, \bTh)] - E[L_Q(\bThh, \bTh)] \over E[L_Q(\X, \bTh)]}, \quad 100\times {E[L_S(\S/n, \bSi)] - E[L_S(\bSih, \bSi)] \over E[L_S(\S/n, \bSi)]},
\\&100\times{E[D_{KL}(\phi(\Ybt|\X, \S/n) ,\phi(\Ybt|\bTh, \bSi))] - E[D_{KL}(\phi(\Ybt|\bThh, \bSih) ,\phi(\Ybt|\bTh, \bSi))] \over E[D_{KL}(\phi(\Ybt|\X, \S/n) ,\phi(\Ybt|\bTh, \bSi))] }.
\end{align*}
The values of PRIR for $\si_{kl}=0.9+0.1\de_{kl}$ and $\si_{kl}=0.5^{|k-l|}$ are reported in Tables \ref{sim1} and \ref{sim2}, respectively. 

\smallskip
Table \ref{sim1} treats the case that the data are strongly correlated, while the data have relatively weak correlations in Table \ref{sim2}.
The performances of the estimators are similar in both tables. 
The values of PRIR for the estimation of the mean matrix $\bTh$ are given in the left columns of the tables. 
When the norm of the mean matrix is small, the improvements of the four estimators are very high.
Comparing the two cases of $q=1$ and $q=1/2$, the improvements are higher in the case of $q=1$, which implies that the shrinkage estimators are more improved when eigenvalues of the mean matrix are more dispersed.
Comparing EB, G1 and GB, we can see that G1 and GB are better for $p>m$ while EB is better for $m>p$.
When the norm of the mean matrix is large, GB performs better than the others in most cases.
The values of PRIR for EB and G1 are always positive, and this fact supports that they are minimax in theory.
The performance of G2 is better than EB and G1 in many cases, but the values of PRIR are negative when $(p,n,m)=(10,10,5), (10,10,25)$.
In fact, the sufficient condition for minimaxity of G2 is not satisfied in these cases.
Since GB is the generalized Bayes estimator and its performance is comparable with EM, the generalized Bayes estimator GB is recommendable. 

\smallskip
The simulation results in estimation of the covariance $\bSi$ are given in the middle columns of the tables.
The performances of G2 and GB are comparable.
When $p>m$, G2 dominates the unbiased estimator in theory, and this result is supported by the simulation. 
When $m>p$, G2 is not guaranteed to dominate the unbiased estimator.
In fact, the value of PRIR of G2 is negative for $(p,n,m)=(10,25,30)$.
All the values of PRIR of G1 are positive and the performance is not bad in both the cases.
The parameters of $\al$ and $\be$ in G1 are different from those in G2, and this suggests that the estimators of the form (\ref{GSi2}) can dominate the unbiased estimator in theory when $m\geq p$.

\smallskip
The simulation results of the simultaneous estimation of $(\bTh, \bSi)$ under the Kullback-Leibler divergence are given in the right columns of the tables, where EM denotes the estimator $(\bThh^{EM}, \S/n)$.
Although G2 and GB have negative values in PRIR for $(p,n,m)=(10,25,30)$ in estimation of $\bSi$, their values in PRIR are positive in estimation of $(\bTh, \bSi)$, because they can borrow the risk gains in the estimation of $\bTh$.
The parameters in GB are chosen from Corollaries \ref{cor:KLminGB1} and \ref{cor:KLminGB}, and GB performs better than the others when the norm of mean matrix is not small.
When the norm of mean matrix is small, the performances of the four estimators are comparable.
Thus, the generalized Bayes estimator GB is recommendable.

\begin{table}[H]
\caption{\small Estimated PRIR($\%$) of the four estimators GB, G1, G2 and EM where singular values of mean are $s_j=s_0+s_0\times (j-1)/(p/5-1)$ for $j\in\{1,\ldots, p/5\}$ and the others are $s_j=s_0/10^{q}$ and covariance is $\si_{kl}=0.9+0.1\de_{kl}$ (EM in estimate of $(\bTh,\bSi)$ denotes the estimator $(\bThh^{EM}, \S/n)$.}
\centering
$
{\renewcommand\arraystretch{1.1}
\begin{array}{c
              c|
              c
              c@{\hspace{1mm}}|
              c@{\hspace{1mm}}
              c@{\hspace{1mm}}
              c@{\hspace{1mm}}
              c|
              c@{\hspace{1mm}}
              c@{\hspace{1mm}}
              c|
              c@{\hspace{1mm}}
              c@{\hspace{2mm}}
              c@{\hspace{2mm}}
              c@{\hspace{2mm}}
              c
             }
s_0&q&\text{$(p,n)$}&\text{$m$} &\multicolumn{4}{c}{\text{estimate of $\bTh$}}&\multicolumn{3}{c}{\text{estimate of $\bSi$}}&\multicolumn{4}{c}{\text{estimate of $(\bTh, \bSi)$}}\\
&&&&\text{EM}&\text{G1}&\text{G2}&\text{GB}&\text{G1}&\text{G2}&\text{GB}&\text{EM}&\text{G1}&\text{G2}&\text{GB}\\
\hline
0&1&\text{$(10,10)$}&\text{$5$}
&17.9 &21.5 &30.7 &20.7 &6.03 &11.3 &10.0  &8.41 &13.3  &20.4  &15.1
\\
&&&\text{$25$}
&26.6 &19.0 &32.0 &24.8 &6.91  &7.61  &7.76 &12.5 &12.6  &19.0  &15.7
\\
&&\text{$(10,25)$}&\text{$5$}
&30.5 &37.2 &47.5 &25.8 &4.33  &5.76  &5.23 &24.3 &30.4  &38.9  &21.6
\\
&&&\text{$20$}
&31.3 &17.7 &25.7 &27.0 &4.38  &3.04  &2.59 &24.8 &14.9  &21.0  &22.0
\\
&&&\text{$30$}
&43.9 &28.1 &37.5 &37.3 &4.81 &-0.387 &-0.197 &34.9 &23.4  &29.7 &29.6
\\
\hline
1&1&\text{$(10,10)$}&\text{$5$}
&16.2 &19.2 &26.2 &19.2 &6.00 &11.1 &10.1  &7.62 &12.2  &18.2  &14.4
\\
&&&\text{$25$}
&24.2 &17.7 &29.8 &23.1 &6.84  &7.86 &7.80 &11.3 &11.9  &18.1  &15.0
\\
&&\text{$(10,25)$}&\text{$5$}
&27.4 &33.4 &42.2 &24.5 &4.27  &5.64  &5.46 &21.8 &27.4  &34.7  &20.6
\\
&&&\text{$20$}
&28.2 &16.5 &24.0 &25.2 &4.74  &3.78  &3.41 &22.4 &14.1  &19.8  &20.7
\\
&&&\text{$30$}
&40.5 &26.6 &35.4 &35.2 &5.45  &0.878  &1.05 &32.1 &22.3  &28.3  &28.2
\\
\hline
10&1&\text{$(10,10)$}&\text{$5$}
&6.26  &7.39  &4.68 &10.2 &5.90 &10.2 &10.3  &2.94  &6.60  &7.59  &10.2
\\
&&&\text{$25$}
&15.1 &12.6 &20.7 &16.3 &6.25  &8.23  &7.47  &7.08  &9.23 &14.1  &11.6
\\
&&\text{$(10,25)$}&\text{$5$}
&10.9 &12.7 &14.0 &14.1 &4.03  &5.03  &6.59  &8.69 &10.9 &12.2 &12.5
\\
&&&\text{$20$}
&16.2 &11.3 &16.2 &17.0 &5.66  &6.12  &6.03 &12.9 &10.1 &14.1  &14.7
\\
&&&\text{$30$}
&26.8 &19.3 &25.5 &25.4 &7.28  &5.22  &5.32 &21.3 &16.9 &21.4  &21.2
\\
\hline
20&1&\text{$(10,10)$}&\text{$5$}
&1.93  &2.08 &-1.65  &3.55 &5.80 &9.31 &9.86  &0.909  &4.04  &4.15  &6.89
\\
&&&\text{$25$}
&8.22  &7.90 &12.2 &10.0 &5.21 &7.92 &6.58  &3.85  &6.47  &9.94  &8.19
\\
&&\text{$(10,25)$}&\text{$5$}
&3.19  &3.45  &2.84  &6.36 &3.90 &4.65 &6.93 &2.53  &3.54  &3.21  &6.48
\\
&&&\text{$20$}
&7.93  &6.98  &9.84 &10.3 &5.75 &7.26 &7.41  &6.30  &6.73  &9.31  &9.67
\\
&&&\text{$30$}
&15.5 &13.4 &17.4 &17.3 &8.23 &8.33 &8.35 &12.3 &12.3 &15.5 &15.5
\\
\hline
1&1/2&\text{$(10,10)$}&\text{$5$}
&15.1 &17.9 &23.6 &18.4 &6.00 &11.0 &10.1  &7.09 &11.6  &16.9 &14.0
\\
&&&\text{$25$}
&23.4 &17.4 &29.3 &22.8 &6.84  &7.92  &7.81 &11.0 &11.8  &17.9  &14.8
\\
&&\text{$(10,25)$}&\text{$5$}
&25.5 &31.0 &38.9 &23.8 &4.25  &5.59  &5.57 &20.2 &25.5  &32.0  &20.1
\\
&&&\text{$20$}
&27.1 &16.2 &23.6 &24.7 &4.83  &3.97  &3.62 &21.5 &13.9  &19.5  &20.4
\\
&&&\text{$30$}
&39.5 &26.3 &35.0 &34.8 &5.60  &1.17  &1.34 &31.40 &22.0  &28.1  &27.9
\\
\hline
10&1/2&\text{$(10,10)$}&\text{$5$}
&0.799 &0.824 &-1.27  &1.10 &5.75 &8.96 &9.26 &0.375 &3.43  &4.13  &5.41
\\
&&&\text{$25$}
&5.12 &4.91  &6.80  &6.02 &3.98 &6.62 &5.21 &2.40 &4.41  &6.71  &5.58
\\
&&\text{$(10,25)$}&\text{$5$}
&1.33 &1.38  &0.980  &2.45 &3.87 &4.54 &6.45 &1.06 &1.89  &1.71  &3.27
\\
&&&\text{$20$}
&4.88 &4.18  &5.66  &5.85 &4.90 &6.66 &6.90 &3.88 &4.33 &5.86  &6.07
\\
&&&\text{$30$}
&9.68 &8.70 &11.0 &10.9 &7.70 &9.04 &9.02 &7.69 &8.49 &10.6 &10.5
\\
\hline
20&1/2&\text{$(10,10)$}&\text{$5$}
&0.199 &0.199 &-0.408 &0.0996 &5.75 &8.76 &8.56 &0.0844 &3.13  &4.44  &4.57
\\
&&&\text{$25$}
&2.74 &2.45  &2.94 &2.88 &2.30 &3.97 &3.05 &1.28 &2.37  &3.49  &2.97
\\
&&\text{$(10,25)$}&\text{$5$}
&0.339 &0.341  &0.229 &0.299 &3.86 &4.47 &5.16 &0.266 &1.06  &1.10  &1.29
\\
&&&\text{$20$}
&2.47 &1.85  &2.38 &2.44 &2.88 &4.03 &4.20 &1.96 &2.07  &2.72  &2.80
\\
&&&\text{$30$}
&4.91 &3.99  &4.75 &4.73 &4.74 &5.98 &5.95 &3.91 &4.14  &5.00  &4.98
\\
\hline
\end{array}
}
$
\label{sim1}
\end{table}

\begin{table}[H]
\caption{\small Estimated PRIR($\%$) of the four estimators GB, G1, G2 and EM where singular values of mean are $s_j=s_0+s_0\times (j-1)/(p/5-1)$ for $j\in\{1,\ldots, p/5\}$ and the others are $s_j=s_0/10^{q}$ and covariance is $\si_{kl}=0.5^{|k-l|}$ (EM in estimate of $(\bTh,\bSi)$ denotes the estimator $(\bThh^{EM}, \S/n)$.}
\centering
$
{\renewcommand\arraystretch{1.1}
\begin{array}{c
              c|
              c
              c@{\hspace{2mm}}|
              c@{\hspace{2mm}}
              c@{\hspace{1mm}}
              c@{\hspace{1mm}}
              c|
              c@{\hspace{1mm}}
              c@{\hspace{1mm}}
              c|
              c@{\hspace{1mm}}
              c@{\hspace{2mm}}
              c@{\hspace{2mm}}
              c@{\hspace{2mm}}
              c
             }
s_0&q&\text{$(p,n)$}&\text{$m$} &\multicolumn{4}{c}{\text{estimate of $\bTh$}}&\multicolumn{3}{c}{\text{estimate of $\bSi$}}&\multicolumn{4}{c}{\text{estimate of $(\bTh, \bSi)$}}\\
&&&&\text{EM}&\text{G1}&\text{G2}&\text{GB}&\text{G1}&\text{G2}&\text{GB}&\text{EM}&\text{G1}&\text{G2}&\text{GB}\\
\hline
0&1&\text{$(10,10)$}&\text{$5$}
&17.9 &21.5 &30.7 &20.7 &6.03 &11.3 &10.0  &8.41 &13.3  &20.4  &15.1
\\
&&&\text{$25$}
&26.6 &19.0 &32.0 &24.8 &6.91  &7.61  &7.76 &12.5 &12.6  &19.0  &15.7
\\
&&\text{$(10,25)$}&\text{$5$}
&30.5 &37.2 &47.5 &25.8 &4.33  &5.76  &5.23 &24.3 &30.4  &38.9  &21.6
\\
&&&\text{$20$}
&31.3 &17.7 &25.7 &27.0 &4.38  &3.04  &2.59 &24.8 &14.9  &21.0  &22.0
\\
&&&\text{$30$}
&43.9 &28.1 &37.5 &37.3 &4.81 &-0.387 &-0.197 &34.9 &23.4  &29.7 &29.6
\\
\hline
1&1&\text{$(10,10)$}&\text{$5$}
&17.7 &21.0 &29.7 &20.4 &6.04 &11.3 &10.1  &8.34 &13.1  &19.9  &14.9
\\
&&&\text{$25$}
&25.9 &18.7 &31.4 &24.4 &6.90  &7.69  &7.78 &12.1 &12.4  &18.8  &15.6
\\
&&\text{$(10,25)$}&\text{$5$}
&29.9 &36.2 &46.1 &25.5 &4.32  &5.73  &5.29 &23.8 &29.6  &37.8  &21.4
\\
&&&\text{$20$}
&30.2 &17.3 &25.2 &26.5 &4.50  &3.27  &2.85 &24.0 &14.7  &20.7  &21.6
\\
&&&\text{$30$}
&42.9 &27.8 &37.0 &36.8 &4.98 &-0.0619  &0.124 &34.1 &23.1  &29.4  &29.2
\\
\hline
10&1&\text{$(10,10)$}&\text{$5$}
&12.6 &14.8 &17.9 &15.6 &5.98 &10.8 &10.2  &5.92 &10.1  &14.1  &12.7
\\
&&&\text{$25$}
&20.4 &15.3 &25.5 &19.9 &6.60  &8.03  &7.67  &9.54 &10.7  &16.2  &13.4
\\
&&\text{$(10,25)$}&\text{$5$}
&21.5 &25.5 &31.4 &19.4 &4.29  &5.58  &6.14 &17.0 &21.1  &26.1  &16.7
\\
&&&\text{$20$}
&23.5 &13.8 &20.0 &21.0 &5.28  &5.03  &4.81 &18.7 &12.1  &17.0  &17.7
\\
&&&\text{$30$}
&34.4 &22.9 &30.3 &30.1 &6.38  &3.07  &3.20 &27.4 &19.5  &24.7  &24.6
\\
\hline
20&1&\text{$(10,10)$}&\text{$5$}
&8.02  &9.41  &8.03 &11.8 &5.91 &10.3 &10.3  &3.77  &7.56  &9.25  &11.0
\\
&&&\text{$25$}
&16.9 &13.5 &22.3 &17.5 &6.39  &8.17  &7.56  &7.92  &9.73 &14.8  &12.2
\\
&&\text{$(10,25)$}&\text{$5$}
&14.1 &16.5 &19.0 &15.7 &4.19  &5.34  &6.56 &11.2 &14.0 &16.2  &13.8
\\
&&&\text{$20$}
&18.7 &12.1 &17.5 &18.3 &5.58  &5.81  &5.68 &14.9 &10.8 &15.1  &15.7
\\
&&&\text{$30$}
&29.1 &20.3 &26.9 &26.7 &6.99 &4.56 &4.67 &23.1 &17.6 &22.3  &22.2
\\
\hline
1&1/2&\text{$(10,10)$}&\text{$5$}
&17.5 &20.8 &29.2 &20.3 &6.05 &11.3 &10.1  &8.21 &13.0  &19.7  &14.9
\\
&&&\text{$25$}
&25.7 &18.6 &31.3 &24.3 &6.90  &7.70  &7.79 &12.1 &12.4  &18.8  &15.5
\\
&&\text{$(10,25)$}&\text{$5$}
&29.6 &35.8 &45.6 &25.4 &4.32  &5.73  &5.31 &23.5 &29.3  &37.4  &21.3
\\
&&&\text{$20$}
&30.0 &17.3 &25.2 &26.4 &4.52  &3.31  &2.89 &23.9 &14.6  &20.7  &21.6
\\
&&&\text{$30$}
&42.8 &27.7 &36.9 &36.7 &5.01 &-0.00774  &0.178 &34.0 &23.1  &29.3  &29.2
\\
\hline
10&1/2&\text{$(10,10)$}&\text{$5$}
&4.45  &4.97  &0.713  &7.78 &5.83 &9.78 &10.2  &2.09  &5.42  &5.51  &9.08
\\
&&&\text{$25$}
&13.1 &11.6 &18.9 &15.0 &6.15 &8.40  &7.45  &6.15  &8.70 &13.3 &11.0
\\
&&\text{$(10,25)$}&\text{$5$}
&7.75  &8.85  &9.00 &12.2 &4.08 &5.06  &7.04  &6.16  &7.87  &8.19 &11.1
\\
&&&\text{$20$}
&13.5 &10.5 &15.1 &15.8 &5.91 &6.72  &6.71 &10.7  &9.55 &13.4 &13.9
\\
&&&\text{$30$}
&23.7 &18.5 &24.4 &24.2 &7.80 &6.31  &6.39 &18.8 &16.3 &20.7 &20.6
\\
\hline
20&1/2&\text{$(10,10)$}&\text{$5$}
&1.11  &1.15 &-1.63  &1.76 &5.77 &9.07 &9.51 &0.516  &3.59  &4.03  &5.86
\\
&&&\text{$25$}
&5.74  &6.07  &8.92  &7.57 &4.63 &7.50 &6.00 &2.69  &5.30  &8.17  &6.74
\\
&&\text{$(10,25)$}&\text{$5$}
&1.98  &2.09  &1.50  &4.04 &3.94 &4.67 &6.90 &1.57  &2.47  &2.15  &4.63
\\
&&&\text{$20$}
&5.03  &5.39  &7.45  &7.74 &5.53 &7.35 &7.57 &3.99  &5.42  &7.43  &7.71
\\
&&&\text{$30$}
&11.4 &10.8 &13.9 &13.8 &8.22 &9.16 &9.16 &9.10 &10.3 &12.9 &12.8
\\
\hline
\end{array}
}
$
\label{sim2}
\end{table}

\section{Concluding Remarks}
\label{sec:remark}
In this paper, we have derived the generalized Bayes shrinkage estimators with closed forms in estimation of the mean and covariance matrices of the multivariate normal distribution.
To establish the dominance results, we have employed the approach of the unbiased risk estimation.
Using this approach, we have obtained the condition for the generalized Bayes estimator to be minimax in estimation of the mean matrix.
In the case of $p>m$, this approach has been also applied to derive the condition under which the generalized Bayes estimator dominated the unbiased estimator in estimation of the covariance matrix.
In the case of $m\geq p$, however, the approach cannot be applied to get any dominance property of the generalized Bayes estimator.
In the framework of simultaneous estimation of the mean and covariance matrices, we can borrow the risk gain in estimation of the mean matrix to establish the dominance property of the generalized Bayes estimators of the mean and covariance matrices, where the simultaneous risk is measured by the Kullback-Leibler divergence.

\medskip
In the case of $m\geq p$, we cannot show the dominance result using the approach of the unbiased risk estimation.
When $p=1$, Maruyama and Strawderman (2006) proved that the generalized Bayes estimator dominates the unbiased estimator.
Their method is based on the direct calculation of the risk function used by Strawderman (1974).
Our simulation results support the dominance property of the generalized Bayes estimator for $m\geq p$.
This suggests that the dominance result may be established by the direct calculation of risk instead of using the unbiased risk estimation.
However, it may be hard, because one needs to calculate some moments related to noncentral multivariate F distributions.

\bigskip
\noindent
{\bf Acknowledgments.}\ \

Research of this work was supported in part by Grant-in-Aid for JSPS Fellows  (19J22203 and 18K11188).

\renewcommand{\thesection}{\Alph{section}}
\setcounter{section}{0}

\end{document}